\documentclass{article}
\usepackage{amsmath}
\usepackage{color}
\usepackage{graphicx} 
\title{Holonomic Gradient Descent for the Fisher--Bingham Distribution on a $d$-dimensional Sphere}
\author{Tamio Koyama\footnote{Department of Mathematics, Kobe University
                     and JST CREST Hibi team} ,
        Hiromasa Nakayama${}^*$,  \\
        Kenta Nishiyama\footnote{School of Management and Information,
	               University of Shizuoka
                       and JST CREST Hibi team} ,
        Nobuki Takayama${}^*$}
\date{October 12, 2012, Revised June 19, 2013}
\def\comment#1{ }
\def\commentt#1{ }
\def\dd{d}
\def\rr{{\rm rank}(I)}
\def\Del{\partial}
\def\pd#1{ \partial_{#1} }
\def\QED{ Q.E.D. \par \bigbreak} 
\newtheorem{theorem}{Theorem}
\newtheorem{proposition}{Proposition}
\newtheorem{lemma}{Lemma}
\newtheorem{example}{Example}
\newtheorem{algorithm}{Algorithm}
\newtheorem{remark}{Remark}

\begin{document}
\maketitle

{\it Abstract}\/.
We propose an accelerated version 
of the holonomic gradient descent 
and apply it
to calculating the maximum likelihood estimate (MLE)
of the Fisher--Bingham distribution on a $\dd$-dimensional
sphere.
We derive a Pfaffian system (an integrable connection) 
and a series expansion 
associated with the normalizing constant with an error estimation.
These enable us to solve some MLE problems up to dimension $\dd=7$
to a specified accuracy.
\bigbreak


\section{Introduction}
Let $x=(x_{ij})$ 
and $y=(y_i)$ be a matrix parameter of size $(\dd+1) \times (\dd+1)$
such that $x_{ij} = x_{ji}$ for $i\not= j$ 
and a vector parameter of size $\dd+1$, respectively. 
We are interested in the Fisher--Bingham probability distribution
\begin{equation*}
\mu(t;x,y,r)|dt| :=  \frac{1}{Z(x, y, r)} \exp \left(  \sum_{1\leq i \leq j \leq \dd+1}x_{ij}t_it_j  +\sum_{i=1}^{\dd+1}y_it_i  \right)|dt|
\end{equation*}
on the $\dd$-dimensional sphere 
$S^\dd(r) = \{ (t_1, \ldots, t_{\dd+1}) \,|\, \sum_{i=1}^{\dd+1} t_i^2 = r^2, r>0 \}$
and the maximum likelihood estimate (MLE) of the parameters $x$ and $y$ of this probability distribution. 
Here, the function $Z$ is the normalizing constant defined as
\begin{equation}  \label{Zdef}
Z(x, y, r) = \int_{S^\dd(r)} \exp \left(  \sum_{1\leq i \leq j \leq \dd+1}x_{ij}t_it_j  +\sum_{i=1}^{\dd+1}y_it_i  \right) |dt|
\end{equation}
and $|dt|$ denotes the standard measure on the sphere of radius $r$ such that
$\int_{S^d(r)} |dt| = r^d\frac{2\pi^{(d+1)/2}}{\Gamma((d+1)/2)}$.

Solving the MLE problem involves finding the maximum of the function in $x, y$
$$ \prod_{k=1}^N \mu(T^{(k)}; x,y,1) $$
for given data vectors $T^{(k)}$, $k=1, \ldots, N$ in the $t$-space 
$S^d(1)$. 
In order to compute the MLE,
we need approximate values for the normalizing constant $Z$
and its derivatives.
In the case of $\dd = 2$, the normalizing constant is expressed in terms of
the Bessel function and 
there are several approaches for computing MLEs
in directional statistics
\cite{kent82}, \cite{MJ}, \cite{wood}.
However, there are few studies 
on approximating the normalizing constant
for the case of $\dd > 2$
and applications to the MLE.
Among these,
Kume and Wood \cite{kume-wood} proposed a method to evaluate the normalizing constant
by utilizing the Laplace approximation of the integral for $\dd > 2$ 
and
Kume and Walker \cite{kume-walker} gave a series approximation of the normalizing constant.

In this paper, we propose a different method for evaluating the normalizing constant
and present applications to the MLE.
Our method is based on the holonomic gradient descent (HGD) proposed in \cite{n3ost2},
which utilizes a holonomic system of linear differential equations satisfied 
by the normalizing constant
and gives the MLE accurately.
The HGD consists of four steps.
The first step is to derive a holonomic system of linear partial differential
equations for the normalizing constant.
The second step is to translate the holonomic system into a Pfaffian system,
which is, roughly speaking, a set of ordinary differential equations 
with respect to the parameters $x_{ij}$ and $y_i$
for the normalizing constant.
These two steps can be performed by a {\it symbolic} computation (the Gr\"obner basis method)
if the size of the problem
is moderate.
The remaining steps utilize {\it numerical} computation.
The third step is to evaluate the normalizing constant and its derivatives
at an initial point.
To do this, we can use numerical integration for a rough evaluation
or a series expansion for a more accurate evaluation.
The last step is to extend the evaluated values to other points
needed for the MLE
by using the Pfaffian system and a numerical solver of ordinary differential equations.

It is shown in \cite{koyama1} and \cite{n3ost2}
that the normalizing constant of the Fisher--Bingham
distribution is a holonomic function 
in $x, y, r$ and consequently that it is annihilated by 
a holonomic ideal for which an explicit expression is given.
This is the first step when applying the HGD.
For the second step,
we need to translate the ideal into a Pfaffian system.
This is performed on a computer for $\dd \leq 2$ in \cite{n3ost2}; however,
this is not possible for $\dd > 2$ on current computers using Gr\"obner basis algorithms
due to the high computational complexity.

In this paper, 
we overcome the difficulty of high complexity and complete the remaining steps
for a general dimension:
we give an accelerated version of the HGD as a general method,
we derive the Pfaffian system of the Fisher--Bingham distribution 
for general $\dd$ 
{\it by hand},
and we derive a series expansion of the normalizing constant with an error estimation.
We will demonstrate that the accelerated version of HGD on our Pfaffian system 
works well up to $\dd=7$
to a specified accuracy for a certain class of problems.
We also propose a general method for evaluating the numerical errors
and apply it to the Fisher--Bingham distribution.

\section{Holonomic Gradient Descent with Pfaffian System of Factored Form}
\label{sec:ff}
The HGD introduced in \cite{n3ost2} 
is a general algorithm for solving MLE problems
for holonomic unnormalized distributions accurately. 
We herein propose an accelerated version of the HGD.
The modification is small, but it allows a drastic performance improvement,
as we will see in the case of the Fisher--Bingham distribution.

In this section, we maintain a general setting to explain
our accelerated method.
A function $f(y_1, \ldots, y_n)$ is called a {\it holonomic function}
when it satisfies an ordinary differential equation with polynomial coefficients
for each variable $y_i$.
In other words, the function $f$ is a holonomic function when
the function is annihilated by an ordinary differential operator of the form
$$ \sum_{j=0}^{r_i} a_{ij}(y) \left(\frac{\partial}{\partial y_i}\right)^j, \quad
  a_{ij}(y) \in {\bf C}[y_1, \ldots, y_n].
$$
By virtue of this ordinary differential operator, the function $f$
can be regarded as a solution of a Pfaffian system discussed below.
An important property of holonomic functions is that the integral
$\int_R f(y_1, \ldots, y_n) dy_n$
is a holonomic function with respect to $y_1, \ldots, y_{n-1}$
under a suitable condition on the integration domain $R$.
In the landmark paper \cite{zeilberger},
Zeilberger introduced the notion of holonomic functions and applied it
to mechanically prove special function identities with systems of linear
partial differential equations and this important property.
It has been proved that the normalizing constant $Z$ of the Fisher--Bingham
distribution is a holonomic function in $x, y, r$ \cite{n3ost2}

Let $f(t;\theta)$ be a holonomic unnormalized probability distribution
with respect to $t$ and $\theta$,
where $\theta=(\theta_1, \ldots, \theta_m)$ is a parameter vector,
and let $Z(\theta) = \int_U f(t;\theta)dt$ be the normalizing constant.
Let $I$ be a holonomic ideal in the ring of differential operators in $\theta$
that annihilates $Z$.
The operator $\partial/\partial \theta_i$ is denoted by $\pd{\theta_i}$.
In order to apply the HGD to MLE problems,
we need an explicit expression for the Pfaffian system
associated with the holonomic ideal $I$
as an input to a numerical solver.
The Pfaffian system is used to numerically evaluate the likelihood function
and its gradient or its Hessian.
Let us review the definition of the Pfaffian system (see, e.g., \cite{n3ost2} for details).
Let $\rr$ be the holonomic rank of the ideal $I$ 
and let $F$ be a vector of the standard monomials of a Gr\"obner basis of $I$.
The length of this vector is $\rr$.
We denote the elements of $F$ by $\pd{}^\alpha$, where $\alpha \in S$.
We assume that the first element of $F$ is $\pd{}^0=1$.
The Pfaffian system is a set of differential operators
which annihilate the vector-valued function
$F(Z) = ( \pd{}^\alpha Z \,|\, \alpha \in S)^T$
that are of the form
$ \pd{\theta_i}-P_i $,
where $P_i$ are $\rr \times \rr$ matrices with rational function entries
which satisfy
$$ \nabla \circ \nabla = 0, \quad 
  \nabla = d - \sum P_i d\theta_i .
$$
In this section, the symbol $d$ in the definition of $\nabla$ 
indicates 
the exterior derivative with respect to the variables $\theta$.

In some of the literature, the definition of the Pfaffian system does not
include the integrability condition   $\nabla \circ \nabla = 0$,
but we will call the integrable Pfaffian system of equations simply
the Pfaffian system for short.
Our Pfaffian system can be regarded as the integrable connection $\nabla$. 

The Pfaffian system can be obtained 
by an algorithmic method explained, e.g.,  in \cite{n3ost2}.
However, it requires heavy computation.
For example,
the computation for the Fisher--Bingham distribution
could be performed within a reasonable time using current computer technology only up to the case of a $2$-dimensional sphere.
Moreover, the Pfaffian system obtained with this method requires heavy
numerical computation in the HGD,
because in general the entries of $P_i$ are huge rational functions.
The latter drawback is removed by using our accelerated version of HGD introduced below.

\begin{algorithm}
\begin{enumerate}
\item Construct a Pfaffian system of the form
\begin{equation} \label{eq:pfaffianQR}
   \pd{\theta_i} - R_i^{-1}(\theta) Q_i(\theta)
\end{equation}
where $Q_i$ and $R_i$ are $\rr \times \rr$ matrices with
{\it polynomial} function entries.
\item Evaluate the normalizing constant $F(Z)$ at an initial parameter $\theta^0$.
\item $k=0$
\item Evaluate the gradient of the likelihood function by using $F(Z)$ at $\theta^k$
     (see \cite{n3ost2}).  
     If the gradient is $0$, then stop.
     Determine the value of the new parameter $\theta^{k+1}$ by standard procedures
    of gradient descent.
    $\theta^{k+1}$ must be sufficiently close to $\theta^k$.
\item Evaluate the approximate value of $F(Z)$ at $\theta^{k+1}$.
      It is, for instance, approximately equal to 
\begin{equation} \label{eq:eulerScheme}
      F(Z)(\theta^k) + 
        \sum_{i=1}^d R_i(\theta^k)^{-1} Q_i(\theta^k) F(Z)(\theta^k) 
          \cdot (\theta^{k+1}_i-\theta^k_i).
\end{equation}
\item $k \rightarrow k+1$ and go to 4.
\end{enumerate}
\end{algorithm}
We call the Pfaffian system of the form (\ref{eq:pfaffianQR})
the Pfaffian system {\it of factored form}.
The main difference between the HGD in \cite{n3ost2} and our proposed method
is (\ref{eq:eulerScheme}).
In the original version, the factored matrix $R_i^{-1} Q_i$ is expressed as a single matrix 
with entries of
(huge) rational functions, but in our proposed method, we express it
as the two matrices $Q_i$ and $R_i$ and
the inverse of $R_i$ is calculated numerically 
in each iteration step.
We note that the approximation in (\ref{eq:eulerScheme}) should be replaced with
a more accurate and efficient numerical scheme such as the Runge--Kutta method.

\begin{remark} \rm
\item When we have a Gr\"obner basis of $I$, the Pfaffian system of the form
(\ref{eq:pfaffianQR}) can be obtained by the computation of normal forms
by the Gr\"obner basis
and then solving linear equations in the ring of polynomials.
This procedure is general, but it can require a huge amount of computational resources.
When we apply this method to problems, 
we need to find shortcuts based on the individual problems in order to solve the problems efficiently.
We will do this for the Fisher--Bingham distribution in the next section.
\end{remark}
\begin{remark} \rm
The matrix $R_i^{-1} Q_i$ may have the form
\begin{equation} \label{eq:pfaffianQR2}
\sum_j R_{ij}^{-1} Q_{ij} T_{ij}^{-1} S_{ij} \cdots
\end{equation}
where
$R_{ij}$, $Q_{ij}$, $S_{ij}$, $T_{ij}$, $\ldots$ 
are $\rr \times \rr$ matrices in polynomial function
entries.
This form is referred to as the multi-factored form
or simply the factored form.
\end{remark}
In the following, we will sometimes call the accelerated version of HGD simply the HGD.

\section{Pfaffian System for the Normalizing Constant}

It is shown in \cite{n3ost2} and \cite{koyama1} 
that the normalizing constant $Z$ in (\ref{Zdef})
of the Fisher--Bingham
distribution is a holonomic function 
in $x, y, r$ and consequently that it is annihilated 
by a holonomic ideal $I$.
The holonomic ideal $I$ is generated by
the following operators in the ring of differential operators.
\begin{eqnarray}
\label{diffopA}
&&\Del_{ij} - \Del_i\Del_j \quad (1 \leq i \leq j \leq \dd + 1), \\ 
\label{diffopB}
&&\sum _{i=1}^{\dd+1} \Del_i^2 - r^2, \\ 
\label{diffopC}
&&x_{ij}\Del_i^2+2(x_{jj}-x_{ii})\Del _i\Del_j
	-x_{ij}\Del_j^2  +\sum _{1 \leq k \leq \dd + 1, k \neq i, j}
	\left( x_{kj}\Del_i\Del_k-x_{ik}\Del_j\Del_k\right) \nonumber \\
&& \quad
        +y_j\Del_i -y_i\Del_j  \quad (1 \leq i < j \leq \dd + 1), \\
\label{diffopE}
&& r\Del_r -2\sum _{1 \leq i\leq j \leq \dd + 1}x_{ij}\Del_i\Del_j
	- \sum_{i=1}^{\dd+1} y_i\Del_i -\dd 
\label{diffopR}
\end{eqnarray}
Here, $\pd{ij}$, $\pd{i}$ 
and $\pd{r}$ stand for $\frac{\partial}{\partial x_{ij}}$, $\frac{\partial}{\partial y_i}$, and
$\frac{\partial}{\partial r}$, respectively.
Note that we assume $x_{ij} = x_{ji}$.

We want to translate these into a Pfaffian system
of the form (\ref{eq:pfaffianQR}) or (\ref{eq:pfaffianQR2})
which is used in the accelerated HGD explained in the previous section.

Before proceeding to the discussion of the general $\dd$-dimensional case, 
we illustrate our method in the case of $\dd=1$ and $r=1$.
Let $I_1$ be the left ideal generated by 
\begin{align}
& \pd{11} - \pd{1}^2,  \quad \pd{12} - \pd{1} \pd{2}, \quad \pd{22} - \pd{2}^2, \\
& \pd{1}^2 + \pd{2}^2 - 1, \label{diffopB1}\\
& x_{12} \pd{1}^2 + 2(x_{22} - x_{11}) \pd{1} \pd{2} - x_{12} \pd{2}^2 +y_2 \pd{1} - y_1 \pd{2} \label{diffopC12} 
\end{align}
in the ring of differential operators.
The holonomic rank of $I_1$ is $4$.
Let $F$ be a vector of operators  $(1,\pd{1},\pd{2}, \pd{1}^2)^T$.
We want to find a matrix $P$ whose entries are rational functions such that 
$\pd{1} F \equiv P F$ holds modulo the left ideal $I_1$.
Here, $ s \equiv t $ means that each element of $s-t$ belongs to $I_1$. 
Since $\pd{1} F = (\pd{1},\pd{1}^2,\pd{1}\pd{2}, \pd{1}^3)^T$,
we need to express $\pd{1}\pd{2}$ and $\pd{1}^3$ in terms of $F$ modulo $I_1$.
Eliminating $\pd{2}^2$ from (\ref{diffopC12}) by (\ref{diffopB1}),
we obtain
\begin{align} \label{dy1dy2}
2(x_{22}-x_{11}) \pd{1}\pd{2} 
&\equiv -x_{12} \pd{1}^2 + x_{12} \underline{\pd{2}^2} - y_2 \pd{1}+y_1 \pd{2} \nonumber \\
&\equiv -x_{12} \pd{1}^2 + x_{12}(1-\pd{1}^2) - y_2 \pd{1}+y_1 \pd{2} 
  \nonumber \\
&\quad\quad (\text{the underlined term is reduced by } (\ref{diffopB1})) \nonumber \\
&= (x_{12}, -y_2, y_1, -2x_{12}) F 
\end{align}
Thus, we have expressed $\pd{1}\pd{2}$ in terms of $F$.
We now try to express $\pd{1}^3$ in terms of $F$.
From $\pd{1} \times (\ref{diffopC12})$,  
we obtain
\begin{align*}
&x_{12} \pd{1}^3 + 2(x_{22} - x_{11}) \pd{1}^2 \pd{2} 
 \equiv x_{12} \pd{1} \underline{\pd{2}^2} - y_2 \pd{1}^2 + y_1 \underline{\pd{1}\pd{2}} + \pd{2}  \\
&\equiv x_{12} \pd{1} (1-\pd{1}^2) - y_2 \pd{1}^2 + \frac{y_1}{2(x_{22} - x_{11})} (x_{12} - y_2 \pd{1} + y_1 \pd{2} - 2 x_{12} \pd{1}^2) + \pd{2} \\
& \quad  \quad (\text{the underlined terms are reduced by } (\ref{diffopB1}) \text{ and } (\ref{dy1dy2})), 
\end{align*}
and consequently we have
$2x_{12} \pd{1}^3 + 2(x_{22} - x_{11}) \pd{1}^2 \pd{2} \equiv (a, b, c, d) F$,
where
\begin{align*} 
&a = \frac{y_1 x_{12}}{2(x_{22} - x_{11})}, \quad
b = x_{12} - \frac{y_1 y_2}{2(x_{22}-x_{11})},  \\
&c = 1 + \frac{y_1^2}{2(x_{22}-x_{11})}, \quad
d = -y_2 - \frac{x_{12}y_1}{x_{22} - x_{11}}.
\end{align*}
By a similar computation for $\pd{2} \times (\ref{diffopC12})$, 
we have 
$$
2(x_{22} - x_{11}) \pd{1}^3 - 2 x_{12} \pd{1}^2 \pd{2} \equiv (a', b', c', d') F, $$
where 
\begin{align*}
&a' = -y_1 + \frac{x_{12}y_2}{2(x_{22}-x_{11})}, \quad
b' = 1 + 2 (x_{22}-x_{11}) - \frac{y_2^2}{2(x_{22}-x_{11})}, \\
&c' = -x_{12} + \frac{y_1 y_2}{2(x_{22}-x_{11})},  \quad 
d' = y_1 - \frac{x_{12} y_2}{x_{22}-x_{11}}.
\end{align*}
Therefore, we have 
$$
\begin{pmatrix}
2 x_{12} & 2(x_{22} - x_{11}) \\
2(x_{22} - x_{11}) & -2 x_{12} 
\end{pmatrix}
\begin{pmatrix}
\pd{1}^3 \\
\pd{1}^2 \pd{2} 
\end{pmatrix}
\equiv 
\begin{pmatrix}
a & b & c & d \\
a' & b' & c' & d'
\end{pmatrix}
F.
$$
Multiplying the both sides by the inverse matrix 
$
\begin{pmatrix}
2 x_{12} & 2(x_{22} - x_{11}) \\
2(x_{22} - x_{11}) & -2 x_{12} 
\end{pmatrix}^{-1}
$,
we can express $\pd{1}^3$ in terms of $F$. 
Thus we have obtained a factored form of $P$ in $\pd{1} F \equiv P F$.
The identity $\pd{1} F \equiv P F$ gives a Pfaffian equation
for the direction $y_1$.
In other words, the differential equation
$$ \frac{\partial F(Z)}{\partial y_1} = P F(Z), \quad 
  F(Z) = \left( Z, \frac{\partial Z}{\partial y_1}, 
          \frac{\partial Z}{\partial y_2}, \frac{\partial^2 Z}{\partial y_1^2}\right)^T
$$ 
holds.
This is an ordinary differential equation 
for the vector-valued function $F(Z)$ with respect to the variable $y_1$.
It is easy to see that $P$ is of the form (\ref{eq:pfaffianQR2}).
Ordinary differential equations for the other directions 
$\pd{2}, \pd{11}, \pd{12}, \pd{22}$
can be obtained analogously.

For the general $\dd$,
let $F$ be the vector of operators
\begin{equation} \label{vecF}
(1, \pd{1}, \ldots, \pd{\dd+1}, \pd{1}^2, \ldots, \pd{\dd}^2)^T.
\end{equation}
\begin{theorem}  \label{H_i}
There exists a $(2d+2) \times (2d+2)$ matrix $H_i$
which has a factored form (\ref{eq:pfaffianQR2})
and satisfies the relation
$\Del_iF \equiv H_iF \ {\rm mod}\, I$.
\end{theorem}

An expression of $H_i$ as a factored form and
a proof of this theorem, which is technical, 
will be given in the Appendix.

The relation between $\pd{ij}F$ and $F$ can be easily obtained by Theorem 
\ref{H_i}.
In fact, since $\pd{ij} F \equiv \pd{i}\pd{j} F \equiv \pd{j}\pd{i}F$
by (\ref{diffopA}),
we have
\begin{equation}
\Del_{ij}F \equiv \Del_j\Del_i F 
\equiv \Del_j (H_i F)
\equiv \frac{\Del H_i}{\Del y_j}F + H_i(\Del_j F)
\equiv \left(\frac{\Del H_i}{\Del y_j} + H_iH_j\right) F.
\end{equation}
We denote by $H_{ij}$ the matrix $\frac{\Del H_i}{\Del y_j} + H_iH_j$.
The matrix such that $\pd{r} F \equiv H_r F$ can be obtained easily
by utilizing (\ref{diffopR}).
Thus, we have obtained the relations
\begin{equation} \label{eq:allPfaffians}
\pd{i} F \equiv H_i F,  \quad 
\pd{ij} F \equiv H_{ij} F, \quad 
\pd{r} F \equiv H_r F.
\end{equation}

In \cite{kn2t2}, we prove that the holonomic rank of
$I$ is equal to $2\dd+2$.
Therefore, the Pfaffian equations are expressed in terms of 
$(2\dd+2) \times (2\dd+2)$ matrices.
Our matrices in (\ref{eq:allPfaffians}) are exactly these matrices. 
The integrability conditions of Pfaffian equations imply
$ 
\frac{\partial H_i}{\partial y_j} + H_i H_j =
\frac{\partial H_j}{\partial y_i} + H_j H_i $.

In \cite{n3ost2}, the differential equations satisfied by the likelihood function
for $\dd=1$ and $\dd=2$ are
derived by a heavy Gr\"obner basis computation and we could not obtain
them for $\dd \geq 3$.
It is known that the Gr\"obner basis computation has the double-exponential complexity with respect to the number of variables
(see, e.g., \cite{bayer-stillman})
and we usually have to avoid deriving Gr\"obner bases by computer
for large problems. Instead, we can sometimes derive
Gr\"obner bases by hand and apply them to interesting applications.
By virtue of Theorem \ref{H_i} for the general dimension, 
we can describe the differential equation satisfied by the likelihood
function with matrices in factored form 
of which factors are relatively small matrices with polynomial entries.
If we calculate the inverse matrices in the factored forms by symbolic computation, 
we would obtain the same result with the Gr\"obner basis method.
In order to apply for the HGD, 
we do not need to calculate these inverse matrices with polynomial entries symbolically;
instead, we need only calculate the inverse matrices numerically when variables are
restricted to real number values in each step of the Runge--Kutta method.
This will become a key ingredient of our algorithm,
which will be discussed in section \ref{algorithm_and_numerical_result}.

\begin{remark} \rm
The matrices $H_i$, $H_{ij}$, $H_r$ have simple forms when $x$ is a diagonal
matrix.
In \cite{so3}, the MLE of the Fisher distribution on $SO(3)$
is obtained by the HGD with differential equations for the normalizing
constant with diagonalized arguments.
It is a natural question to ask whether a simplification analogous to the diagonal $x$
case is possible.
Unfortunately, an analog of Lemma 2 of \cite{so3} does not hold
except for the case of $y=0$.
It is possible to evaluate the gradient of $Z$ 
from values of $Z$ for diagonal $x$
by Proposition \ref{prop:diag} 
given later;
however, this requires computation of a transformation matrix 
to diagonalize the matrix $x$ at each step of the gradient descent.
On the other hand, we do not need to do this computation for the diagonal form
in the HGD by the Pfaffian system for the full parameters $x, y, r$.
\end{remark}

\section{Series Expansion for the Normalizing Constant}

Let us define the function $\tilde{Z}$ by the integral
\begin{equation} \label{tilde-z}
\tilde{Z}(\tilde{x},\tilde{y},\tilde{r}) = 
 \int_{S^\dd(r)} \exp \left( \sum_{i=1}^{\dd+1}(\tilde{x_i}t_i^2+\tilde{y_i}t_i) \right) |dt|.
\end{equation}
The function satisfies the invariance relation
\begin{equation} \label{tilde-z-invariance}
 \tilde{Z}(\tilde x,\tilde y,1) = \tilde{Z}(r^{-2} \tilde x, r^{-1} \tilde y, r).
\end{equation}
This function is  the restriction of the normalizing constant
$Z$ to the diagonalized ${\tilde x}$.
Since the normalizing constant is invariant under the action of the orthogonal 
group $O(\dd+1)$,
we can express $F(Z)$  in terms of $F(\tilde{Z})$.
The following proposition can be obtained by a straightforward calculation.

\begin{proposition}  \label{prop:diag}
Suppose that the real symmetric matrix $x$
is diagonalized by an orthogonal matrix
$P=(p_{ij})$, and
set
$\tilde{x}= P^TxP,\, \tilde{y} = P^Ty,\, \tilde{r} = r$.
Then, we have
\begin{eqnarray*}
Z(x,y,r) &=& \tilde{Z}(\tilde{x}, \tilde{y}, \tilde{r})
\\
\frac{\Del Z}{\Del y_i}(x,y,r) &=& \sum_{k=1}^{\dd+1} p_{ik}\frac{\Del \tilde{Z}}{\Del \tilde{y}_k}(\tilde{x},\tilde{y},\tilde{r})
\\
\frac{\Del^2 Z}{\Del y_i^2}(x,y,r) &=& 
\sum_{k=1}^{\dd+1} p_{ik}^2 \frac{\Del^2 \tilde{Z}}{\Del \tilde{y}_k^2}(\tilde{x},\tilde{y},\tilde{r})
\\
&&
-\sum_{1\leq k<\ell\leq \dd+1} \frac{p_{ik}p_{i\ell}}{\tilde{x}_k- \tilde{x}_\ell}
    \left(
	     \tilde{y}_k \frac{\Del \tilde{Z}}{\Del \tilde{y}_\ell   }(\tilde{x},\tilde{y},\tilde{r})
	    -\tilde{y}_\ell\frac{\Del \tilde{Z}}{\Del \tilde{y}_k}(\tilde{x},\tilde{y},\tilde{r})
	\right).
\end{eqnarray*}
\end{proposition}

Kume and Walker give a series approximation of
the Fisher--Bingham distribution and consequently that of the normalizing
constant ${\tilde Z}({\tilde x},{\tilde y},1)$ \cite{kume-walker}.
Their expression is easily rescaled to the case including the
parameter
$r$.
We will give an error estimate of this series approximation.
We will omit the tilde symbol (`\~{}') for $x$ and $y$ in the following 
where this should cause no confusion.
\begin{theorem} \label{series}
\begin{enumerate}
\item {\rm \cite{kume-walker}}
The restricted normalizing constant has the following series expansion:
\begin{equation}  \label{eq:series}
\tilde{Z}(x,y,r)
=
S_d\cdot
\sum_{\alpha, \beta \in {\bf N}_0^{d+1}}
r^{d+2|\alpha+\beta |}
\frac{(d-1)!!\prod_{i=1}^{d+1}(2\alpha_i+2\beta_i-1)!!}{(d-1+2|\alpha|+2|\beta|)!!\alpha !(2\beta)!}
x^\alpha y^{2\beta}.
\end{equation}
Here, $S_d = \int_{S^d(1)}|dt|$ denotes the surface area of the $d$-sphere of radius $1$,
and ${\bf N}_0 = \{0, 1, 2, \dots \}$.
For a multi-index $\alpha \in {\bf N}_0^{d+1}$, we define 
$\alpha! = \prod_{i=1}^{d+1}\alpha_i!,\, \alpha!! = \prod_{i=1}^{d+1} \alpha_i!!,$
and $|\alpha| = \sum_{i=1}^{d+1} \alpha_i.$
\item The truncation error of the series is estimated as
\begin{eqnarray}
&&
  \left|
  S_d \cdot
  \sum_{|\alpha+\beta| \geq N}
  r^{d+2|\alpha+\beta |}
  \frac{(d-1)!!\prod_{i=1}^{d+1}(2\alpha_i+2\beta_i-1)!!}{(d-1+2|\alpha|+2|\beta|)!!\alpha !(2\beta)!}
  x^\alpha y^{2\beta} 
  \right|  \nonumber
\\
&\leq&
S_d \cdot
\frac{r^d}{N!} 
\left(
r^2\sum_i(|x_i|+|y_i|^2)
\right)^N 
\frac{N+1}{N+1-r^2\sum_i(|x_i|+|y_i|^2)}
\end{eqnarray}
when $N$ is sufficiently large.
\end{enumerate}
\end{theorem}

We note that the series (\ref{eq:series}) converges slowly
when $r^2\sum_i(|x_i|+|y_i|^2)> 1$
and converges relatively rapidly
when $r^2\sum_i(|x_i|+|y_i|^2) \leq 1$.
The derivatives of $\tilde Z$ are expressed as derivatives of
the right-hand side of (\ref{eq:series}).

{\it Proof of 2}\/.
We have the following estimates
\begin{eqnarray*}
&&
  \left|
  \sum_{|\alpha+\beta| \geq N}
  r^{d+2|\alpha+\beta |}
  \frac{(d-1)!!\prod_{i=1}^{d+1}(2\alpha_i+2\beta_i-1)!!}{(d-1+2|\alpha|+2|\beta|)!!\alpha !(2\beta)!}
  x^\alpha y^{2\beta} 
  \right|
\\
&\leq& 
  \sum_{|\alpha+\beta| \geq N}
  \left|
  r^{d+2|\alpha+\beta |}
  \frac{(d-1)!!\prod_{i=1}^{d+1}(2\alpha_i+2\beta_i-1)!!}{(d-1+2|\alpha|+2|\beta|)!!\alpha !(2\beta)!}
  x^\alpha y^{2\beta} 
  \right|
\\
&\leq &
  \sum_{|\alpha+\beta| \geq N}
  \left|
  r^{d+2|\alpha+\beta |}
  \frac{1}{\alpha !(2\beta)!} x^\alpha y^{2\beta} 
  \right|
\\
&\leq &
  \sum_{|\alpha+\beta| \geq N}
  \left|
  r^{d+2|\alpha+\beta |}
  \frac{1}{\alpha !\beta!} x^\alpha y^{2\beta} 
  \right|
\leq 
  r^d
  \sum_{n \geq N}
  \frac{r^{2n}}{n!}
  \sum_{|\alpha+\beta|=n}
  \frac{n!}{\alpha !\beta!} |x|^\alpha |y|^{2\beta} 
\\
&\leq& 
  r^d
  \sum_{n \geq N}
  \frac{r^{2n}}{n!}
  \left(
  \sum_{i=1}^{d+1} (|x_i|+|y_i^2|)
  \right)^n
\leq
  r^d
  \sum_{n \geq N}
  \frac{1}{n!}
  \left(
  r^2\sum_{i=1}^{d+1} (|x_i|+|y_i^2|)
  \right)^n.
\end{eqnarray*}
Set $L(x,y,r) =   r^2\sum_{i=1}^{d+1} (|x_i|+|y_i^2|)$.
Assume that $N$ is sufficiently large
so that $L/(N+1) < 1$.
We have the estimate
\begin{equation}
  r^d
  \sum_{n \geq N}
  \frac{1}{n!}
  L(x,y,r)^n
\leq
\frac{r^d}{N!} 
L(x,y,r)^N \frac{N+1}{N+1-L(x,y,r)}
\end{equation}
 by the estimate
\begin{eqnarray*}
\sum_{n\geq N}  \frac{1}{n!}L^n
&=&
\frac{1}{N!} L^N \sum_{n\geq N}  \frac{1}{(N+1)_{n-N}}L^{n-N} 
\\
&=&
\frac{1}{N!}L^N \sum_{n=0}^{\infty}  \frac{1}{(N+1)_n}L^n 
\\
&\leq&
\frac{1}{N!}L^N \sum_{n=0}^{\infty}  \frac{1}{(N+1)^n}L^n 
\\
&=&
\frac{1}{N!}L^N \frac{N+1}{N+1-L}.
\end{eqnarray*}
\QED

\section{Numerical Evaluation of the Normalizing Constant}

In order to efficiently evaluate $\tilde{Z}$ numerically,
we use 
the holonomic gradient method (HGM)
(see, e.g., \cite{hnt2}).
The HGM is a method for evaluating the normalizing constant 
by utilizing a system of differential equations.
In the case of the Fisher--Bingham distribution, 
we numerically evaluate the series (\ref{series})
in the domain $r^2\sum_i(|x_i|+|y_i|^2) \leq 1$
and extend the numerical evaluation outside this domain by using
a differential equation with respect to $r$.
To use this method, we prepare the following theorem.

\begin{theorem}
\begin{enumerate}
\item
The function $\tilde{Z}$ is annihilated by the left ideal $\tilde{I}$
generated by
\begin{eqnarray*}
 &&A_i=\Del_{y_i}^2-\Del_{x_{i}}  \quad  (1\leq i \leq \dd+1), \\
 &&B    =\Del_{y_1}^2+\cdots+\Del_{y_{\dd+1}}^2-r^2, \\
 &&C_{ij}=2(x_{i}-x_{j})\Del_{y_i}\Del_{y_j}+y_i\Del_{y_j}-y_j\Del_{y_i}  \quad (1\leq i < j \leq \dd+1),\\ 
 &&E    =r\Del_r-2\sum_{i=1}^{\dd+1}x_{i}\Del_{y_{i}}^2-\sum_{i=1}^{\dd+1}y_i\Del_{y_i}-\dd.
\end{eqnarray*}
\item 
Set $\tilde{F} = (\Del_{y_1}, \dots, \Del_{y_{\dd+1}}, \Del_{y_1}^2, \dots, \Del_{y_{\dd+1}}^2)^T$.
Then, we have
$ \pd{r} \tilde{F} \equiv P^{(r)} \tilde{F} \ {\rm mod}\, \tilde{I}$.
Here, the matrix
$P^{(r)}= (p_{ij}^{(r)})$ is defined by 
\begin{eqnarray*}
rp_{ij}^{(r)}       &=&                     (2x_ir^2+1)\delta_{ij}      + \sum_{k=1}^{\dd+1} y_i\delta_{j(k+\dd+1)} \quad (1 \leq i \leq \dd+1 ),\\
rp_{(i+\dd+1)j}^{(r)} &=& y_ir^2\delta_{ij} + (2x_ir^2+2)\delta_{j(i+\dd+1)}+ \sum_{k\neq i}      \delta_{j(k+\dd+1)} \quad (1 \leq i \leq \dd+1 )
\end{eqnarray*}
for $1 \leq j \leq 2\dd+2$.
\end{enumerate}
\end{theorem}

The proof of this theorem is analogous to that for the non-diagonal $x$ case.
\comment{{\it Add a proof later.}}

\begin{example} \rm
In the case of $\dd=1$, the matrix $P^{(r)}$ is
$$
\frac{1}{r} 
\left(
\begin{array}{cccc}
2r^2x_1+1 &         0 &       y_1 &       y_1 \\
        0 & 2r^2x_2+1 &       y_2 &       y_2 \\
   r^2y_1 &         0 & 2r^2x_1+2 &         1 \\
        0 &    r^2y_2 &         1 & 2r^2x_2+2 \\  
\end{array}
\right).
$$
\end{example}

We note that the largest eigenvalue of $P^{(r)}$ is $O(r)$.
Our implementation of the HGM numerically solves the ordinary differential equation
\begin{equation} \label{eq:hgm-r}
 \frac{\partial G}{\partial r} = (P^{(r)} - r \lambda E) G
\end{equation}
instead of solving 
$\pd{r} {\tilde F}({\tilde Z}) = P^{(r)} {\tilde F}({\tilde Z})$,
where
$\lambda$ is the largest eigenvalue of $\lim_{r \rightarrow +\infty} P^{(r)}/r$
and the vector-valued function $G$ is defined by
$\tilde F({\tilde Z}) = \exp(\lambda r^2/2) G$.
This scalar scaling is necessary, because the adaptive Runge--Kutta method
requires an absolute error bound of the solution 
to automatically make meshes finer,
and when a solution grows exponentially,
meshes become too small to maintain an absolute error bound.

Now let us discuss the accuracy of the HGM.
The truncation error of the series approximation is estimated in Theorem
\ref{series}.
We want to estimate the numerical error 
caused by applying the HGM.
In other words, we want to estimate how much the truncation error is
magnified 
by solving the ordinary differential equation numerically.
We propose a practical method to do this.
Note that this method can be applied to any HGM, but we will explain it 
in the case of the Fisher--Bingham distribution.
In this method, we assume that initial values are governed 
by a probability measure.
This assumption is natural in, e.g., molecular modeling, and
some classes of non-linear ordinary differential equations are  
studied under this assumption (see, e.g., \cite{cssw}).
In our case, the ordinary differential equation is linear 
and the problem is much easier.
Since we have not found a reference relevant to our case,
we include below a discussion on the behavior of solutions of a linear
equation under random initial data.
For a given initial value vector $Z_0$ at $r=r_0$,
we denote by $R Z_0$ 
the output obtained at $r=r_1$ by solving the ordinary differential equation 
where we may suppose that $R$ is a constant matrix,
because the Runge--Kutta solver can be regarded as a linear map from the
input to the output
under the assumption that round-off errors and cancellation errors by
floating point arithmetic are sufficiently small
and that the automatic mesh refinement process is fixed.
When $Z_0$ is regarded as a random
vector distributed as a multivariate normal distribution,
the output $R Z_0$ is also a random vector distributed as a multivariate
 normal distribution. 
In our implementation, 
we perturb $Z_0$ with random numbers of which the standard deviation
is $\varepsilon/2$, where $\varepsilon$ is a truncation error,
 and solve the ordinary differential equation for
these perturbed initial values.
We evaluate the mean and the standard deviation 
of the first component of $R Z_0$
and these give an evaluation of the normalizing constant and its
statistical error bound.
This bound is more practical than that by interval arithmetic.

\begin{figure}[tbp]
 \begin{center}
  \includegraphics[width=100mm]{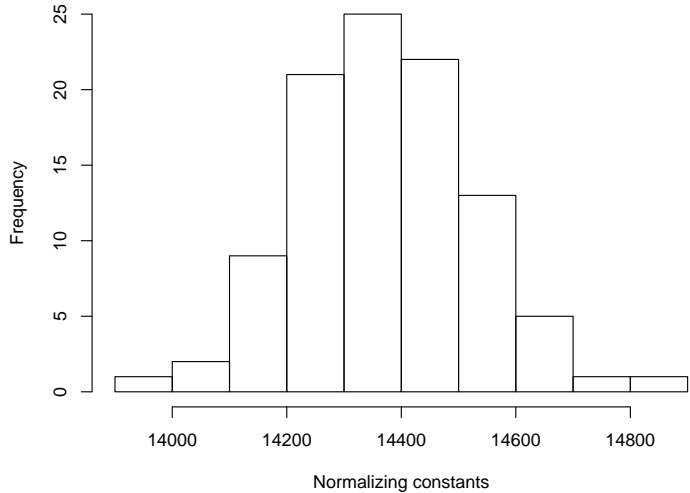}
 \end{center}
 \caption{Histogram of normalizing constants by the HGM with random initial values}
 \label{fig:hist1}
\end{figure}

For example, Figure \ref{fig:hist1} shows a histogram of the
normalizing constant which is generated by the above procedure with 
$\varepsilon/2 = 0.1$,
$r_0 = 1$, $r_1= \sum |x_i| + \sum y_i^2$,
$\dd=3$, $x={\rm diag}(1.2,2.5,3.2,3.6)$,
$y=(2.3,5.3,4.2,0.1)$.
The series is evaluated at $x/r_1^2,y/r_1,r=1$ and extended to $r=r_1$. 
In this case, the standard derivation of the normalizing constant is
evaluated as $156.6288$ and the confidence interval with probability
$0.95$ is $[   14065.6,    14679.6]$

Kume and Wood \cite{kume-wood} gave a saddle point approximation of the
normalizing constant for the Fisher--Bingham distribution.
Our method evaluates the normalizing constant with an error bound.
Table \ref{fig:2} shows values by the HGM and by the third-order saddle point approximation of
Kume and Wood.
Here, $\dd=4$ and 
$(x_{ij})={\rm
diag}(x_{11},2x_{11},3x_{11},4x_{11},5x_{11})$,
$ 0.5 \leq x_{11} \leq 10 $, 
$(y_k)=(0.5 y_0,
0.4y_0,0.3y_0,0.2y_0,0.1y_0)$, $y_0=3$.
The absolute error bound to solve 
(\ref{eq:hgm-r}) by the adaptive Runge--Kutta method is set to 
$10^{-6} \sum_{i=1}^{2\dd+2} G_i/(2\dd+2)$ 
and $\varepsilon$ is $10^{-5}$.
The values of the standard deviation imply that the values by the HGM
have at least 6-digit accuracy with 95\% confidence.

\begin{table}[htb]
\begin{center}
\caption{Normalizing constants}  \label{fig:2}
\begin{tabular}{|r|r|r||r|} \hline
 $x_{11}$ &\multicolumn{2}{|c||}{HGM} & \multicolumn{1}{c|}{\makebox[82pt][c]{Kume-Wood, 3}}\\ \cline{2-3}
             & \makebox[82pt][r]{HGM} & \makebox[82pt][r]{standard deviation}    &    \\ \hline \hline
 0.5  &         189.243&	 1.737976e-04& \underline{189}.763                 \\ \hline 
 1.0  &         985.529&	 9.102497e-04& \underline{9}94.043                 \\ \hline 
 1.5  &         5856.78&	 5.424156e-03& \underline{58}08.16                 \\ \hline 
 2.0  &         39075.8&	 3.624707e-02& \underline{3}7602.6                 \\ \hline 
 2.5  &          287231&	 2.667160e-01& \underline{2}71557                  \\ \hline 
 3.0  &      	2.28420e+06&	 2.122623e+00& \underline{2}.15158e+06             \\ \hline 
 3.5  &  	1.93448e+07&	 1.798630e+01& \underline{1}.82924e+07             \\ \hline 
 4.0  &  	1.72236e+08&	 1.602082e+02& \underline{1}.63939e+08             \\ \hline 
 4.5  &  	1.59584e+09&	 1.484901e+03& \underline{1.5}2931e+09             \\ \hline 
 5.0  &  	1.52663e+10&	 1.420891e+04& \underline{1}.4717e+10              \\ \hline 
 5.5  &  	1.49868e+11&	 1.395204e+05& \underline{1.4}5179e+11             \\ \hline 
 6.0  &  	1.50274e+12&	 1.399244e+06& \underline{1}.46123e+12             \\ \hline 
 6.5  &  	1.53345e+13&	 1.428082e+07& \underline{1}.49556e+13             \\ \hline 
 7.0  &  	1.58797e+14&	 1.479060e+08& \underline{1.5}5222e+14             \\ \hline 
 7.5  &  	1.66504e+15&	 1.551038e+09& \underline{1.6}302e+15              \\ \hline 
 8.0  &  	1.76459e+16&	 1.643961e+10& \underline{1.7}299e+16              \\ \hline 
 8.5  &  	1.88748e+17&	 1.758618e+11& \underline{1.8}5223e+17             \\ \hline 
 9.0  &  	2.03531e+18&	 1.896519e+12& 1.99905e+18                          \\ \hline 
 9.5  &      	2.21040e+19&	 2.059834e+13& \underline{2}.1716e+19              \\ \hline 
 10.0 &  	2.41579e+20&	 2.251392e+14& \underline{2}.37462e+20             \\ \hline 
\end{tabular}
\end{center}
\end{table}

\section{Algorithm and Numerical Results}
\label{algorithm_and_numerical_result}

In \cite[Algorithm 1 and Theorem 2]{n3ost2}, we give
an algorithm to obtain the MLE
for the Fisher--Bingham distribution.
This algorithm is valid for general dimensions, but it cannot be used 
for more than two dimensions with the current level of computer technology because of
the high computational complexity of the Gr\"obner basis computation.
We replace the Gr\"obner basis computation part
(steps 1, 2, 3 in \cite[Algorithm 1]{n3ost2}) with
our derivation of the Pfaffian system of factored form given 
in the Theorem \ref{H_i}, and
replace the numerical integration of (\ref{Zdef}) part
with the evaluation by the series (\ref{series})
and extend values to slowly convergent domains of the series
by the HGM.
For efficiency, we calculate the inverse matrices 
in our expressions for $H_{ij}$ and $H_i$ numerically
during the steps of the adaptive Runge-Kutta method
as explained in Section \ref{sec:ff} regarding the accelerated version of the HGD.
This enables us to solve maximum likelihood estimation problems
in more than two dimensions case with the HGD.
More precisely, we have the following complexity result.

\begin{theorem} \label{complexity_theorem}
The complexity of the series expansion method,
the HGM, and the HGD
for the Fisher--Bingham distribution on the $\dd$-dimensional sphere is
\begin{eqnarray*}
&&  O( (2\dd+2)^{N+1}/N! ) + \mbox{\rm (complexity of solving the ODE with respect to $r$)} \\
&&+ 
  O((2\dd+2)^3) \times (\mbox{\rm steps of the convergence of gradient descent}).
\end{eqnarray*}
The first and the second terms are the complexity to evaluate the initial values
$F(Z)$ 
up to degree $N$ and 
the third term is the complexity of the HGD.
\end{theorem}
{\it Proof}\/.
The number of terms of the truncated series of (\ref{eq:series})
is ${{2\dd+2+N} \choose {2\dd+2}}={{2\dd+2+N} \choose {N}} = O((2\dd+2)^N/N!)$.
The coefficients of the series can be evaluated by a recursive relation.
We need $2\dd+1$ derivatives of ${\tilde Z}$.
Thus, we obtain the first term.

Our HGD requires the computation of the inverses of $(2\dd+2) \times (2\dd+2)$ matrices
in each step of the HGD,
by Theorem \ref{H_i}.
This corresponds to the third term.  
\QED

\comment{
Our computational experiments suggest that
giving a practical bound for 
the complexity of numerically solving the ordinary differential equation
with respect to $r$ with a prescribed accuracy
is not easy.
}

We implemented our algorithm firstly in Maple and
next in the C language by using 
the GNU scientific library \cite{gsl}.
The prototype written in Maple is useful for debugging our C code.
Our C code is automatically generated by our code generation
program {\tt pfn\_gen\_c\_2.rr},
which can be obtained from the URL in the Example \ref{ex:3},
on {\tt Risa/Asir} \cite{risa/asir}.

In order to apply the HGD, we need to find a good starting point $\theta^0=(x^0,y^0)$.
We use the following method.
\begin{algorithm} \label{alg:two-steps} \rm
\begin{enumerate}
\item Take a random point ${\tilde \theta}^0=(x,y)$ satisfying
$0<x_{ij}<1$ and $0< y_k <1$.
\item Apply the Nelder--Mead algorithm,
which does not require the gradient and may also be replaced with other methods,
to find an approximate optimal point $\theta^0$ of the likelihood function
from the starting point ${\tilde \theta}^0$
(see, e.g., \cite{nocedal-wright}).
Normalization constants are evaluated by the HGM,
which may be replaced by other methods.
\item Apply the HGD with the starting point $\theta^0$.
If the HGD stops normally, we are done.
If the HGD stops at $\theta^1$ because of a numerical instability, go to step 2
with ${\tilde \theta}^0 = \mbox{ a point in a neighborhood of } \theta^1$.
\end{enumerate} 
\end{algorithm}
An alternative and heuristic way to avoid the retry in the last step
is to abort the computation when a numerical instability occurs
and then restart the algorithm with a new randomly chosen starting point.
This procedure can be implemented in parallel.

We present some examples to illustrate the performance of our new algorithm
and its implementation.

\begin{example} \rm
The problem ``Astronomical data'' given in \cite{n3ost2}
is solved in 2.58 seconds on a 32-bit virtual machine,
the host machine of which is an Intel Xeon E5410 (2.33G Hz)
processor based computer.
In contrast, our implementation in \cite{n3ost2} spends
17.3 seconds for the HGD and more than an hour for computing a Gr\"obner basis
and deriving a Pfaffian system.
\end{example}
The following timing data are taken on the same machine.
Values in the following table are given in seconds. 


\begin{example} \rm \label{ex:3}
Problem names beginning with sk\_'s in Table \ref{fig:table3}
are problems on the $k$-dimensional sphere.
These problems are generated by a random number generator
according to the Fisher--Bingham distribution.
We choose $8$ random points in the parameter space as starting points 
for Algorithm \ref{alg:two-steps}.
The HGD aborts $2$ times in the $8$ tries in the worse case on $S^3$.
This rate increases to $7$ aborts of $8$ tries in the worse case on $S^7$.
The timing data shown are those for the first successful HGD among the $8$ starting points.
The first step time in the table is that of the step of applying the Nelder--Mead algorithm
with the HGM.
\begin{table}[htb]
\begin{center}
\begin{tabular}{|c|c|c|c|}
\hline
Problem   & Time (1st step) & Time of the HGD and steps & Total \\ \hline
\hline
s3\_e1 & 9.8 & 3.2(73) & 13  \\
s3\_e3 & 10  & 3.1(66) & 13  \\
\hline
s4\_e1 & 50  & 14(93)  & 64  \\
s4\_e2 & 50  & 28(183) & 78  \\
s4\_e3 & 50  & 11(75)  & 61  \\
\hline
s5\_e1 & 220 & 80(142) & 300 \\
s5\_e2 & 221 & 140(121)& 361 \\
s5\_e3 & 222 & 66(117) & 288 \\
\hline
s6\_e1 & 828 & 247(172)& 1075 \\
\hline
s7\_e1 & 2679& 571(183)& 3250 \\
\hline
\end{tabular}
\caption{Performance of the accelerated HGD} \label{fig:table3}
\end{center}
\end{table}
These sample data and programs are obtainable from our web page.
\footnote{\tt http://www.math.kobe-u.ac.jp/OpenXM/Math/Fisher-Bingham-2}
\end{example}

\section{Conclusion and Open Problems}

We show that the HGD can solve some MLE problems up to dimension $\dd=7$
by utilizing an explicit expression of the 
the Pfaffian system of factored form
and the series expansion of the normalizing constant.
However, there are two problems in applying our method 
efficiently to arbitrary data.
\begin{enumerate}
\item
In examples, we find a starting point for applying the HGD
by the Nelder-Mead algorithm with the HGM for evaluating the normalization constant.
This method seems to work well for our examples, but it is not very efficient.
Finding a good starting point efficiently for arbitrary data is an open question.
\item There are domains in $(x,y)$-space where the normalizing constant
cannot be evaluated to a given accuracy within a reasonable time by the HGM,
because the normalizing constant is huge in these domains,
which includes domains where $|y|$ is large.
\end{enumerate}
Although, there still remain important open problems,
our proposed method evaluates the normalizing constant
and its derivatives to a specified accuracy for a sufficiently broad set of parameters
and solves MLE problems. 
We can easily control the accuracy of evaluations of the normalizing constant
and so it is possible to apply our method for evaluations to other approximation
methods.

\appendix
\section{Proof to the Theorem \ref{H_i}}

We define two auxiliary vectors of operators to present the expression.
We sort the set of the square free second order operators
$$
\{\Del_i\Del_j | 1 \leq i<j \leq \dd+1 \}
$$
by the lexicographic order.
This gives a vector of operators of the length $\dd(\dd+1)/2$:
\begin{equation}
F^{(2)} = (\Del_1\Del_2,\, \Del_1\Del_3, \dots, \Del_\dd\Del_{\dd+1} )^T.
\end{equation}
We sort the set of the third order operators
$$
\{\Del_i\Del_j\Del_k | 1 \leq i \leq j \leq k \leq \dd+1,\, j \leq \dd \}
$$
by the lexicographic order.
We denote by $F^{(3)}$ the sorted vector
\begin{equation}
F^{(3)} = 
(\Del_1\Del_1\Del_1 ,\,\Del_1\Del_1\Del_2,\dots,\Del_1\Del_1\Del_{\dd+1},\,  \Del_1\Del_2\Del_2, \dots, \Del_\dd\Del_\dd\Del_{\dd+1})^T.
\end{equation}
The length of this vector $\dd(\dd + 1)(\dd + 5)/6$ is denoted by $m$.

When two operators $\ell_1$ and $\ell_2$ are the same modulo the ideal $I$,
we denote it by $\ell_1 \equiv \ell_2$. 
By examining the proof of Lemmas 2 and 3 of \cite{n3ost2}, we obtain the following
two lemmas which give an expression of the second and the third order
operators $F^{(2)}$ and $F^{(3)}$ in terms of $F$.
\begin{lemma} \label{lemma1}
We have
\begin{equation}
P^{(2)} F^{(2)} + Q^{(2)} F \equiv 0.
\end{equation}
Here, $P^{(2)}$ is an invertible $\dd(\dd+1)/2 \times \dd(\dd+1)/2$ matrix
and $Q^{(2)}$ is a $\dd(\dd+1)/2 \times (2\dd+2)$ matrix
of which entries are as follows.
\begin{eqnarray*}
P_{ij,kl}^{(2)}&=&
\begin{cases}
2(x_{jj}-x_{ii}) & (i=k, j=l)\\
x_{jl} & (i=k, j\neq l) \\
x_{jk} & (i= l , j\neq k)\\
-x_{ik} & (i\neq k , j = l)\\
-x_{il} & (i\neq l , j = k)
\end{cases}
\\
Q_{ij,k}^{(2)}&=&
\begin{cases}
y_j\delta_{k,i+1}-y_i\delta_{k,j+1}+x_{ij}\delta_{k,i+\dd+2}-x_{ij}\delta_{k,j+\dd+2} & (j\leq \dd)\\
y_j\delta_{k,i+1}-y_i\delta_{k,j+1}+x_{ij}\delta_{k,i+\dd+2}
-r^2x_{i,\dd+1}\delta_{k1} + \sum_{\ell = 1}^\dd x_{i,\dd+1}\delta_{k,\ell+\dd+2} & (j = \dd+1)
\end{cases}
\end{eqnarray*}
\end{lemma}
Here, $\delta$ is Kronecker's $\delta$ and 
$P_{ij,kl}^{(2)}$ is the matrix element of $P^{(2)}$ standing for
$\partial_i\partial_j$ and $\partial_k \partial_l$ in $F^{(2)}$.
We use this notation of the index of the matrix element in the sequel.

\begin{lemma} \label{lemma2}
We have
\begin{equation}
P^{(3)} F^{(3)} + Q^{(3)} F^{(2)} + R^{(3)}F \equiv 0.
\end{equation}
Here, $P^{(3)}$, $Q^{(3)}$, and $R^{(3)}$ are an invertible
$m \times m$ matrix, an $m \times \dd(\dd+1)/2$
matrix and an $m \times (2\dd+2)$ matrix of polynomial entries respectively.
Entries are defined as follows.
\begin{eqnarray*}
P^{(3)}_{ijk,abc}
&=&
\begin{cases}
(\delta_{k,\dd+1}+1)x_{jk}\delta_{ai}\delta_{bj}\delta_{cj} \\ \quad\quad
+2(x_{kk}-x_{jj})\delta_{ai}\delta_{bj}\delta_{ck}
+(\delta_{k,\dd+1}-1)x_{jk}\delta_{ai}\delta_{bk}\delta_{ck} \\ \quad\quad
+\sum _{l \neq j, k}
\left(
x_{kl}\delta'_{abc; ijl}
-x_{jl}\delta'_{abc; ikl}
+x_{jk}\delta_{k,\dd+1}\delta'_{abc;ill}
\right)
& (i \leq j < k \leq \dd+1)
\\
x_{ij}\delta_{ai}\delta_{bi}\delta_{cj}
+2(x_{jj}-x_{ii})\delta_{ai}\delta_{bj}\delta_{cj}
-x_{ij}\delta_{aj}\delta_{bj}\delta_{cj} \\ \quad\quad
+\sum _{l \neq i, j}
\left(
x_{jl}\delta'_{abc;ijl}  
-x_{il}\delta'_{abc;jjl} 
\right)
& (i < j = k < \dd+1)
\\
\sum _{s=1}^\dd (
x_{s,\dd+1}\delta'_{abc;is,\dd+1}                
-2(x_{\dd+1,\dd+1}-x_{ii}) \delta'_{abc;iss}     
\\ \quad\quad
+\sum _{l \neq i}
x_{il}\delta'_{abc;lss} 
)
& (i = j = k < \dd+1)
\\
\end{cases}
\\
Q^{(3)}_{ijk,ab}
&=&
\begin{cases}
(1-\delta_{ij})y_k\delta_{ai}\delta_{bj} -y_j\delta_{ai}\delta_{bk}  & (i \leq j < k \leq \dd+1)\\
y_j\delta_{ai}\delta_{bj}                           & (i < j = k < \dd+1)\\
y_{\dd+1}\delta_{ai}\delta_{b,\dd+1}                    & (i = j = k < \dd+1)\\
\end{cases}
\\
R^{(3)}_{ijk,a}
&=&
\begin{cases}
-x_{jk}r^2\delta_{k,\dd+1}\delta_{a,i+1}-\delta_{ij}\delta_{a,k+1}+y_k\delta_{ij}\delta_{a,i+\dd+2}& (i \leq j < k \leq \dd+1)\\
-y_i\delta_{a,j+\dd+2} +\delta_{a,i+1}& (i < j = k < \dd+1)\\
-y_ir^2\delta_{a1}
+(2(x_{\dd+1,\dd+1}-x_{ii})r^2+1)\delta_{a,i+1} \\ \quad\quad
-\sum_{l\neq i}x_{il}r^2\delta_{a,l+1}
+\sum_{l<\dd+1}y_i\delta_{a,l+\dd+2}
& (i = j = k < \dd+1)
\end{cases}
\end{eqnarray*}
where 
$$
\delta'_{abc;ijk} = 
\begin{cases}
1 & (\Del_a\Del_b\Del_c = \Del_i\Del_j\Del_k) \\
0 & (\Del_a\Del_b\Del_c \neq \Del_i\Del_j\Del_k) 
\end{cases}.
$$
\end{lemma}

{\it Proof}\/.
We denote by $C_{ij}$ the differential operator (\ref{diffopC}) in $I$;
we put
$$
C_{ij}=x_{ij}\Del_i^2+2(x_{jj}-x_{ii})\Del _i\Del_j
	-x_{ij}\Del_j^2
        +\sum _{k \neq i, j}
	\left( x_{kj}\Del_i\Del_k-x_{ik}\Del_j\Del_k\right)
	+y_j\Del_i -y_i\Del _j.
$$
Define a differential operator $G_{ijk}$ ($i \leq j \leq k \leq \dd+1, j \leq \dd$)
by
$$
G_{ijk}=
\begin{cases}
\Del_iC_{jk} & (i \leq j < k \leq \dd+1),\\
\Del_jC_{ij} & (i < j= k \leq \dd), \\
\Del_{\dd+1}C_{i,\dd+1} & (i = j =  k \leq \dd).
\end{cases}
$$
We expand $G_{ijk}$ in the ring of differential operators and express it
in terms of the elements of $F$, $F^{(2)}$, and $F^{(3)}$.
For example, when $i < j < k < \dd+1$, we have
\begin{eqnarray*}
G_{ijk}&=& \Del_iC_{jk}\\
&=&
\Del_i \big(x_{jk}\Del_j^2+2(x_{kk}-x_{jj})\Del _j\Del_k
	-x_{jk}\Del_k^2 \\
&& \quad
        +\sum _{l \neq j, k}
	\left( x_{lk}\Del_j\Del_l-x_{jl}\Del_k\Del_l\right)
	+y_k\Del_j -y_j\Del _k \big)\\ 
&=&
x_{jk}\Del_i\Del_j^2+2(x_{kk}-x_{jj})\Del_i\Del _j\Del_k
	-x_{jk}\Del_i\Del_k^2 \\
&&\quad
        +\sum _{l \neq j, k}
	\left( x_{kl}\Del_i\Del_j\Del_l-x_{jl}\Del_i\Del_k\Del_l\right)
	+y_k\Del_i\Del_j -y_j\Del_i\Del _k
\end{eqnarray*}
which yields
$P^{(3)}_{ijk,ijj}, 
P^{(3)}_{ijk,ijk}, \ldots,
Q^{(3)}_{ijk,ij}, 
Q^{(3)}_{ijk,ik}$. 
Analogous expansions and rewritings for the other cases
give the conclusion.
\QED

We denote by 
${\rm Mat}(k,l,S)$
the space of 
the $k\times l$ matrices with entries in the set $S$.
Let ${\bf Q}[x,y,r]$ denote the ring of polynomials with coefficients in ${\bf Q}$.
\begin{lemma} \label{pfaffian}
The vector $F$ satisfies the identity
\begin{equation}\label{eqM}
A\Del_iF \equiv BF + CF^{(2)} +EF^{(3)}.
\end{equation}
Here, 
$
A=(a_{pj})\in{\rm Mat}(2\dd+2,2\dd+2,{\bf Q}[x,y,r]), 
B=(b_{pj})\in{\rm Mat}(2\dd+2,2\dd+2,{\bf Q}[x,y,r]), 
C=(c_{p,jk})\in{\rm Mat}(2\dd+2,\dd(\dd+1)/2,{\bf Q}[x,y,r]), 
E=(e_{p,jk\ell})\in{\rm Mat}(2\dd+2,m,{\bf Q}[x,y,r]) 
$
and 
$A$ is invertible in the space of the matrices with entries in 
the field of rational functions ${\bf Q}(x,y,r)$.
Explicit expressions of these matrices are given in 
(\ref{eq:abce1}), (\ref{eq:abce2}), (\ref{eq:abce3}), (\ref{eq:abce4}),
(\ref{eq:abce5}), (\ref{eq:abce6}), (\ref{eq:abce7}), (\ref{eq:abce8}).
Note that $A,B,C,E$ depend on the index $i$.
\end{lemma}
Notation: $c_{p,jk}$ means the element at the $p$-th row of $C$
and the column of $C$ standing for 
$\partial_j \partial_k =\partial_k \partial_j$.
$e_{p,jk\ell}$ is defined analogously.

{\it  Proof}\/.
The both sides of $(\ref{eqM})$
is a column vector of the length $2\dd+2$.
We will determine the rows of $A, B,C, E$ from generators of $I$.
Note that the index $i$ is fixed over the proof.

{\it The first rows}\/.
The first element of the vector $\Del_i F$ is
$\Del_i$, then we have
\begin{equation} \label{eq:abce1}
a_{11} = 1, \quad b_{1,i+1}=1
\end{equation}
and the other elements of the first rows of $A, B, C, E$ are $0$.

{\it The $(j+1)$-th rows $(1 \leq j \leq \dd, i \neq j)$}\/.
Using the differential operator $(\ref{diffopC})$ in $I$,
we have
$$
x_{ij}\Del_i^2
+2(x_{jj}-x_{ii})\Del _i\Del_j
+\sum _{k \neq i, j} x_{kj}\Del_i\Del_k
\equiv
x_{ij}\Del_j^2
+\sum _{k \neq i, j} x_{ik}\Del_j\Del_k
+y_i\Del _j 
-y_j\Del_i.
$$
Therefore, we may put as
\begin{eqnarray} \label{eq:abce2}
&&a_{j+1,i+1}=x_{ij},\quad a_{j+1,j+1}=2(x_{jj}-x_{ii}), \\
&&a_{j+1,k+1}=x_{kj}\quad (1 \leq k \leq \dd+1, k \neq i, k \neq j), \nonumber\\
&&b_{j+1,j+1}=y_i,\quad b_{j+1,i+1}=-y_j,\quad b_{j+1,j+\dd+2}=x_{ij}, \nonumber \\
&&c_{j+1,jk}=x_{ik} \quad (1 \leq k \leq \dd+1, k \neq i, k \neq j). \nonumber
\end{eqnarray}
Notation: when an index is out of bound, ignore the setting.
For example, we set $b_{j+1,j+\dd+2}=x_{ij}$ when $j+d+2 \leq 2d+2$.
The other elements of the $(j+1)$-th rows of $A, B, C, E$ are $0$.

{\it The $(i+1)$-th rows}\/.
The $(i+1)$-th element of the vector $\Del_i F$ is
$\Del_i^2$.
When $i \leq \dd$, we put
\begin{equation} \label{eq:abce3}
a_{i+1,i+1} = 1, \quad b_{i+1,i+\dd+2}=1
\end{equation}
and the other elements of the $(i+1)$-th rows are $0$.
When $i = \dd+1$,
we consider the operator $(\ref{diffopB})$ in the ideal $I$.
Then, we have 
$$
\Del_{\dd+1}^2 \equiv r^2 - \sum _{k=1}^\dd \Del_k^2 
$$
and hence we put
\begin{eqnarray} \label{eq:abce4}
&&a_{\dd+2,\dd+2} = 1, \\
&&b_{\dd+2,1}=r^2, \quad b_{\dd+2,k+\dd+2} = -1 \quad (1 \leq k \leq \dd). \nonumber
\end{eqnarray}
The other elements of the $(i+1)$-th rows of $A, B, C, E$ are $0$.

{\it The $(\dd+2)$-th rows}\/.
When $i = \dd+1$, it is reduced to the case of the $(i+1)$-th rows.
We assume that $i \leq \dd$.
Using the operators $(\ref{diffopC})$ and $(\ref{diffopB})$ in $I$,
we have
\begin{eqnarray*}
&&
x_{i,\dd+1}\Del_i^2
+2(x_{\dd+1,\dd+1}-x_{ii})\Del _i\Del_{\dd+1}
+\sum _{k \neq i, \dd+1} x_{k,\dd+1}\Del_i\Del_k\\
&\equiv&
x_{i,\dd+1}\Del_{\dd+1}^2
+\sum _{k \neq i, \dd+1} x_{ik}\Del_{\dd+1}\Del_k
+y_i\Del_{\dd+1} 
-y_{\dd+1}\Del_i\\
&\equiv&
x_{i,\dd+1}r^2
-\sum_{k=1}^\dd x_{i,\dd+1}\Del_k^2
+\sum_{k \neq i, \dd+1} x_{ik}\Del_{\dd+1}\Del_k
+y_i\Del_{\dd+1} 
-y_{\dd+1}\Del_i.
\end{eqnarray*}
Hence, we put
\begin{eqnarray} \label{eq:abce5}
&&a_{\dd+2,i+1}=x_{i,\dd+1},\quad a_{\dd+2,\dd+2}=2(x_{\dd+1,\dd+1}-x_{ii}),\quad \\
&& a_{\dd+2,k+1}=x_{k,\dd+1}\quad (1 \leq k \leq \dd, k \neq i),\nonumber \\
&&b_{\dd+2,1}=x_{i,\dd+1}r^2,\quad b_{\dd+2,\dd+2}=y_i,\quad b_{\dd+2,i+1}=-y_{\dd+1}, \nonumber \\
&&b_{\dd+2,l+\dd+2}=-x_{i,\dd+1},\quad (1 \leq l \leq \dd), \nonumber \\
&&c_{\dd+2,k(\dd+1)}=x_{ik} \quad (1 \leq k \leq \dd, k \neq i). \nonumber
\end{eqnarray}
The other elements of the $(\dd+2)$-th rows of $A,B,C,E$ are $0$.

{\it The $(j+\dd+2)$-th rows $(1 \leq j \leq \dd, i \neq j)$}\/.
Using the operator $(\ref{diffopC})$
multiplied by $\Del_j$ from the left hand side,
we have
$$
-2(x_{jj}-x_{ii})\Del_i\Del_j^2	
\equiv 
x_{ij}\Del_i^2\Del_j -x_{ij}\Del_j^3
+\sum _{k \neq i, j}
\left( x_{kj}\Del_i\Del_j\Del_k-x_{ik}\Del_j^2\Del_k\right)
+y_j\Del_i\Del_j  -y_i\Del_j^2 +\Del_i.
$$
When $i \leq \dd$, we put
\begin{eqnarray} \label{eq:abce6}
&&a_{j+\dd+2,j+\dd+2}=-2(x_{jj}-x_{ii}),\\
&&b_{j+\dd+2,i+1}=1, \quad b_{j+\dd+2,j+\dd+2}=-y_i, \nonumber \\
&&c_{j+\dd+2, ij} = y_j, \nonumber \\
&&e_{j+\dd+2, iij} = x_{ij}, \quad e_{j+\dd+2, jjj}=-x_{ij},\nonumber \\
&& e_{j+\dd+2, ijk} = x_{kj},
 \quad e_{j+\dd+2, jjk} = -x_{ik}
 \quad (1\leq k \leq \dd+1, k\neq i, k\neq j). \nonumber
\end{eqnarray}
The other elements in the $(j+\dd+2)$-th rows of $A,B,C,E$ are $0$.

When $i = \dd+1$, 
We use the operator $(\ref{diffopB})$ and obtain
\begin{eqnarray*}
&&
-2(x_{jj}-x_{\dd+1,\dd+1})\Del_{\dd+1}\Del_j^2 \\
&\equiv& 
x_{\dd+1,j}r^2\Del_j
-2x_{\dd+1,j}\Del_j^3
+y_j\Del_{\dd+1}\Del_j  
-y_{\dd+1}\Del_j^2 
+\Del_{\dd+1}\\
&+&
\sum _{k \neq \dd+1, j} 
\left( 
x_{kj}\Del_{\dd+1}\Del_j\Del_k-x_{\dd+1,k}\Del_j^2\Del_k -x_{\dd+1,j}\Del_j\Del_k^2 
\right).
\end{eqnarray*}
Therefore, we may put as
\begin{eqnarray}  \label{eq:abce7}
&&a_{j+\dd+2,j+\dd+2}=-2(x_{jj}-x_{\dd+1,\dd+1}),\\
&&b_{j+\dd+2,j+1}=x_{\dd+1,j}r^2, \quad b_{j+\dd+2,\dd+2}=1, \quad b_{j+\dd+2,j+\dd+2}=-y_{\dd+1}, \nonumber \\
&&c_{j+\dd+2, j(\dd+1)} = y_j, \nonumber \\
&&e_{j+\dd+2,jjj} = -2x_{ij},\nonumber \\
&&e_{j+\dd+2,ijk} = x_{kj}, \quad e_{j+\dd+2,jjk} = -x_{ik}, \quad e_{j+\dd+2,jkk} = -x_{ij}
 \quad (1\leq k \leq \dd, k\neq j). \nonumber 
\end{eqnarray}
The other elements of the $(j+\dd+2)$-th rows of $A,B,C,E$ are $0$.

{\it The $(i+\dd+2)$-th rows}\/.
We may assume that $i \leq \dd$.
Since the $(i+\dd+2)$-th element of the vector $\Del_i F$ is $\Del_i^3$,
we put
\begin{equation} \label{eq:abce8}
a_{i+\dd+2,i+\dd+2} = 1, \quad e_{i+\dd+2,iii} = 1.
\end{equation}
The other elements of the $(i+\dd+2)$-th rows of $A,B,C,E$ are $0$.
\QED

From the Lemmas \ref{lemma1}, \ref{lemma2}, \ref{pfaffian}, 
we have the Theorem \ref{H_i}, which gives a differential equation
satisfied by the normalizing constant with respect to the variable $y_i$.
As we remarked in the Lemma \ref{pfaffian}, we note that $A, B, C, E$ depend
on the index $i$ and we omit to denote the dependency.

{\it Proof of the Theorem \ref{H_i}}\/.
\begin{eqnarray*}
\Del_iF 
&\equiv& A^{-1}(BF + CF^{(2)} +EF^{(3)}) 
\quad \mbox{ by the Lemma \ref{pfaffian}}\\
&\equiv& A^{-1}
  \left(
     BF 
    +CF^{(2)} 
    -E(P^{(3)})^{-1} \left(  Q^{(3)}F^{(2)}+R^{(3)}F   \right)
  \right)\\
&& \quad \quad \mbox{by the Lemma \ref{lemma2}} \\
&\equiv& A^{-1}
  \left(
     BF 
    -C(P^{(2)})^{-1}Q^{(2)}F
    -E(P^{(3)})^{-1} \left(  -Q^{(3)}(P^{(2)})^{-1}Q^{(2)}F+R^{(3)}F   \right)
  \right)\\
&&\quad \quad \mbox{by the Lemma \ref{lemma1}} \\
&\equiv& A^{-1}
  \left(
     B 
    -C(P^{(2)})^{-1}Q^{(2)}
    +E(P^{(3)})^{-1}\left(  Q^{(3)}(P^{(2)})^{-1}Q^{(2)}-R^{(3)}   \right)
  \right)F.
\end{eqnarray*}
\QED

\begin{example}\rm
In the case of $\dd=1$ and for the $y_1$ direction,
these matrices are as follows.
\begin{eqnarray*}
F &=& \begin{pmatrix}
 1&  \partial_{{1}}&  \partial_{{2}}&   \partial_{{1}}^{ 2} \\
\end{pmatrix}^T, \\
F^{(2)} &=& \begin{pmatrix}
  \partial_{{1}}  \partial_{{2}}\\
\end{pmatrix}, 
\quad %
F^{(3)} = \begin{pmatrix}
  \partial_{{1}}^{ 3} &    \partial_{{1}}^{ 2}   \partial_{{2}}\\
\end{pmatrix}^T,
\end{eqnarray*}
\begin{eqnarray*}
A &=& \begin{pmatrix}
 1& 0& 0& 0 \\
0&  1& 0& 0 \\
0&  {x}_{12}&   -  2  {x}_{11}+  2  {x}_{22}& 0 \\
0& 0& 0&  1 \\
\end{pmatrix}, 
\quad
B =  \begin{pmatrix}
0&  1& 0& 0 \\
0& 0& 0&  1 \\
   {r}^{ 2}   {x}_{12}&  - {y}_{2}&  {y}_{1}&  - {x}_{12} \\
0& 0& 0& 0 \\
\end{pmatrix},
\\ %
C &=& \begin{pmatrix}
0 \\
0 \\
0 \\
0 \\
\end{pmatrix},
\quad
E  =  \begin{pmatrix}
0& 0 \\
0& 0 \\
0& 0 \\
 1& 0 \\
\end{pmatrix},
\end{eqnarray*}
\commentt{
For  the $y_2$ direction, we have
\begin{eqnarray*}
A &=& \begin{pmatrix}
 1& 0& 0& 0 \\
0&    2  {x}_{11}-  2  {x}_{22}&  {x}_{12}& 0 \\
0& 0&  1& 0 \\
0& 0& 0&   -  2  {x}_{11}+  2  {x}_{22} \\
\end{pmatrix}
\\ %
B &=& \begin{pmatrix}
0& 0&  1& 0 \\
0&  {y}_{2}&  - {y}_{1}&  {x}_{12} \\
  {r}^{ 2} & 0& 0&  - 1 \\
0&    {r}^{ 2}   {x}_{12}&  1&  - {y}_{2} \\
\end{pmatrix}
\\ %
C &=& \begin{pmatrix}
0 \\
0 \\
0 \\
 {y}_{1} \\
\end{pmatrix}
\\ %
E &=& \begin{pmatrix}
0& 0 \\
0& 0 \\
0& 0 \\
 -  2  {x}_{12}& 0 \\
\end{pmatrix}
\end{eqnarray*}
} 
\begin{eqnarray*}
P^{(2)}&=&\begin{pmatrix}
  -  2  {x}_{11}+  2  {x}_{22} \\
\end{pmatrix},
\quad 
Q^{(2)} =\begin{pmatrix}
 -   {r}^{ 2}   {x}_{12}&  {y}_{2}&  - {y}_{1}&   2  {x}_{12} \\
\end{pmatrix},
\end{eqnarray*}
\begin{eqnarray*}
P^{(3)}&=&\begin{pmatrix}
   2  {x}_{11}-  2  {x}_{22}&   2  {x}_{12} \\
  2  {x}_{12}&   -  2  {x}_{11}+  2  {x}_{22} \\
\end{pmatrix},
\quad 
Q^{(3)} = \begin{pmatrix}
 {y}_{2} \\
 - {y}_{1} \\
\end{pmatrix},
\\ %
R^{(3)}&=&\begin{pmatrix}
 -   {r}^{ 2}   {y}_{1}&   -   2   {r}^{ 2}   {x}_{11}+    2  {x}_{22}   {r}^{ 2} + 1&  -   {r}^{ 2}   {x}_{12}&  {y}_{1} \\
0&  -   {r}^{ 2}   {x}_{12}&  - 1&  {y}_{2} \\
\end{pmatrix}.
\end{eqnarray*}
\end{example}

\end{document}